\documentclass{amsart}
\usepackage{graphicx} 
\usepackage{amsmath, amsxtra, amsthm, amsfonts, amssymb, mathtools, mathrsfs, stmaryrd}
\usepackage{hyperref}
\usepackage{url}
\usepackage{tikz}
\usepackage{tikz-cd}
\usetikzlibrary{shapes.geometric,positioning}
\usepackage{array}
\usepackage{enumerate}
\usepackage{cleveref}
\usepackage{quiver}
\usepackage{booktabs}
\usepackage[margin=1in]{geometry}
\usepackage[bbgreekl]{mathbbol}
\usepackage{amsbsy}
\usepackage{wasysym}
\usepackage{hyperref}
\usepackage{bbm}

\DeclareMathAlphabet{\dutchcal}{U}{dutchcal}{m}{n}

\usepackage[style=alphabetic]{biblatex}
\newtheorem{Thm}{Theorem}[section]
\newtheorem{Cor}[Thm]{Corollary}
\newtheorem{Lm}[Thm]{Lemma}
\newtheorem{Prop}[Thm]{Proposition}

\theoremstyle{definition}
\newtheorem{Df}[Thm]{Definition}
\newtheorem*{conj}{Conjecture}
\newtheorem{Ex}[Thm]{Example}

\crefname{Thm}{theorem}{theorems}
\Crefname{Thm}{Theorem}{Theorems}

\crefname{Cor}{corollary}{corollaries}
\Crefname{Cor}{Corollary}{Corollaries}

\crefname{Lm}{lemma}{lemmas}
\Crefname{Lm}{Lemma}{Lemmas}

\crefname{Prop}{proposition}{propositions}
\Crefname{Prop}{Proposition}{Propositions}

\crefname{Df}{definition}{definitions}
\Crefname{Df}{Definition}{Definitions}

\crefname{Ex}{example}{examples}
\Crefname{Ex}{Example}{Examples}

\crefname{Rmk}{remark}{remarks}
\Crefname{Rmk}{Remark}{Remarks}

\crefname{conj}{conjecture}{conjectures}
\Crefname{conj}{Conjecture}{Conjectures}

\theoremstyle{remark}

\newtheorem*{rmk}{Remark}

\title{On the Failure of Unimodality and Top-Heaviness for matroid quotient counts}
\author{Cooper Kofron}
\address{Georgia Institute of Technology}
\email{ckofron3@gatech.edu}
\date{October 2026}

\begin{document}

\maketitle
\begin{abstract}
In this paper, we construct realizable matroids with arbitrarily many peaks in their quotient counts, and matroids for which these counts violate top-heaviness. Failure of top-heaviness in particular resolves a privately communicated conjecture of June Huh. This provides evidence and counterevidence towards a philosophy of \cite{wang2026tropical} relating the quotient poset of a matroid $M$ with its lattice of flats, $\mathcal{L}(M)$.
\end{abstract}
\section{Introduction}

Let $M$ be a matroid of rank $d$ on a ground set $E$. A subset
$F \subseteq E$ is called a \emph{flat} of $M$ if adding any element of $E\setminus F$ to $F$ would increase its rank in $M$. The collection of flats of $M$, denoted by $\mathcal{L}(M)$, forms a lattice under inclusion. We denote $\mathcal{L}_k(M)$ as the set of rank-$k$ flats of $M$, $\mathcal{L}^k(M)$ as the set of corank-$k$ flats of $M$, and the numbers $W_k(M)=|\mathcal{L}_k(M)|$ are the \emph{Whitney numbers of the second kind}. $\mathcal{L}(M)$ is one of the most well-studied objects in matroid theory, and many guiding conjectures have centered around its properties. Two such conjectures were formulated by Rota and Dowling-Wilson in the 1970s:

\begin{conj}[\cite{rota1970combinatorial}] Let $M$ be a matroid. Then $W_k(M)$ forms a unimodal sequence in $k$.
\end{conj} 
\begin{conj}[\cite{dowling1975whitney}] Let $M$ be a matroid of rank $d$. Then
$W_k(M)$ is top-heavy: $W_{k}(M)\leq W_{d-k}(M)$ for $k\leq \frac{d}{2}$.
\end{conj}

Monotonicity for the first half of $W_k(M)$ was shown in the case of realizable matroids in  \cite{huh2017enumeration}, and generally in \cite{braden2020singular}. Unimodality in full was recently disproved in \cite{divouxlarson}. In fact, the failure of unimodality can be arbitrarily bad:

\begin{Thm} [\cite{divouxlarson}] \label{Thm: DLLWPeaks}
    For every positive integer $m$, there is a realizable matroid $M$ whose sequence $W_k(M)$ has exactly $m$ peaks.
\end{Thm}

In contrast, the Dowling-Wilson conjecture was shown to hold by Braden et al.:

\begin{Thm}[\cite{braden2020singular}]
For every matroid $M$ of rank $d$ and every pair of integers
$k \leq j \leq d-k$,
\[
    W_k(M) \leq W_j(M).
\]
\end{Thm}

\emph{Quotients} are a fundamental concept in matroid theory, closely related to that of flats: \begin{Df}[Quotient] \label{Quotient cryptomorphisms}
    Let $M$ and $N$ be matroids on the same ground set. Then $N$ is a \emph{quotient} of $M$ (denoted $M\twoheadrightarrow N$) if $\mathcal{L}(N)\subseteq \mathcal{L}(M)$. Said otherwise, every flat of $N$ is a flat of $M$.
\end{Df} As above, we define $\mathcal{Q}(M)$ as the set of quotients of $M$, $\mathcal{Q}_k(M)$ as the set of rank $k$ quotients of $M$, and $q_k(M)=|\mathcal{Q}_k(M)|$. $\mathcal{Q}(M)$ can also be equipped with a partial order given by $N\leq M\Longleftrightarrow M \twoheadrightarrow N$. Motivated by the Baker--Bowler--Lorscheid theory of \emph{matroids with coefficients}, the recent thesis of Wang \cite{wang2026tropical} introduced the following Main Speculation: $\mathcal{Q}(M)$ should behave similarly to $\mathcal{L}(M)$. He developed evidence both for and against this speculation:

\begin{center}
\begin{tabular}{|c|c|c|}
\hline
 & $\mathcal{L}(M)$ & $\mathcal{Q}(M)$ \\
\hline
Graded & \checkmark & \checkmark \\
\hline
Shellable & \checkmark & \checkmark \\
\hline
Lattice & \checkmark & Not always \\
\hline
Upper semimodular & \checkmark & Not always \\
\hline
\end{tabular}
\end{center}
Importantly, there also is an (order-reversing) injection $\mathcal{L}(M)\hookrightarrow\mathcal{Q}(M)$ given by $F\mapsto (M/F)\oplus U_{0,F}$. In this note, we provide further evidence and counter-evidence towards this Main Speculation by disproving analogs of \emph{both} Rota's and Dowling-Wilson's conjectures for $\mathcal{Q}(M)$. \begin{Thm} \label{Thm: intro_compilation} For each positive integer $m$, there is a realizable matroid $M$ whose sequence of quotient counts has at least $m$ peaks.
\end{Thm}

\begin{Thm}
    For each $d\geq 5$, there exists a matroid $M$ of rank $d$ with $q_2(M)> q_{d-2}(M)$.
\end{Thm}

The counterexample to unimodality is a modification of the example found in \cite{divouxlarson}, supporting Wang's Main Speculation. While it's possible that the matroid from that paper itself has non-unimodal quotient count, the sheer size of their examples makes verifying non-unimodality of $q_k(M)$ computationally an astronomical feat far exceeding current capabilities. The key idea of the counterexample is instead to increase the size of the matroid further (in a minimally invasive way), by taking the direct sum with a large uniform matroid of rank 2. In contrast, our counterexamples to top-heaviness are surprisingly small (the smallest is of rank 5 on 30 elements), and they rely on properties of modular matroids shown in \cite{wang2024lorentzian,wang2026tropical}.

\begin{rmk} \label{Rmk: Pasture_Motivation}
    Since the quotients of $M=U_{n,n}$ are precisely the \emph{labelled} matroids on $[n]$, unimodality of $q_k(M)$ for general $M$ can be thought of as a generalization of the longstanding ``unimodality of matroids" conjecture due to Welsh \cite{PiffWelsh71}. This states that for a fixed ground set $[n]$, the number of rank $k$ matroids on $[n]$ (up to isomorphism) is unimodal in $k$. A different generalization of this conjecture asserted unimodality for the number of rank $k$ $P$-matroids on $[n]$, where $P$ is a fixed finite pasture. That generalization also fails, and the counterexample curiously also involves uniform matroids \cite{KofronLowen26}. It was that disproof, combined with Wang's analogy, which led the author to consider the role of uniform matroids in the current context. 
\end{rmk}
\subsection{Acknowledgements}
The author thanks Shouda Wang, Jui-Hui Chung, and Chayim Lowen for helpful discussions and contributions in the early stages of this project, and for much feedback on the readability of the proof. He also thanks Jayden Wang and Matthew Baker for feedback on a previous draft.

\subsection{Statement on AI use} 
These counterexamples were discovered by AI. The author had observed the importance of the matroid in \cite{divouxlarson}, but quickly realized that computing $q_k(M)$ explicitly was infeasible. He had also observed the importance of uniform matroids in a possible disproof of unimodality, in analogy with \cite{KofronLowen26}. Using AI, Jui-Hui Chung found a violation to unimodality along these lines and communicated it to the author. The present work contains an alternate proof, completed without the use of AI. The original failure of top-heaviness was discovered via LLM by Chayim Lowen, while he was providing feedback on an earlier draft of the current paper. Those counterexamples were likewise broadened, simplified, and strengthened through the author's numerous discussions with Chayim Lowen.

\section{Preliminaries}
For the rest of the paper, we assume all matroids have finite ground sets. We first compile a few elementary results about quotients, starting with cryptomorphisms.

\begin{Prop}\label{Prop: quotient_crypto} \cite{oxley2006matroid} Let $M$ and $Q$ be matroids on the same ground set $E$. Then the following are equivalent: 
    \begin{enumerate}
        \item [1.] $M \twoheadrightarrow Q$.
        \item [2.] $\text{cl}_M(X)\subseteq \text{cl}_Q(X)$ for every subset $X\subseteq E$.
        \item [3.] For every $X\subseteq Y\subseteq E$, $$r_M(Y)-r_M(X)\geq r_Q(Y)-r_Q(X).$$
    \end{enumerate}
    Furthermore, if $M \twoheadrightarrow Q$, then $r_Q(X)\leq r_M(X)$ for every $X\subseteq E$, and $r(M)=r(Q)$ if and only if $M=Q$.
\end{Prop} 

Since there is only one matroid of rank $0$ on any ground set, $q_0(M)=1$ for any matroid $M$. The final observation implies that $q_{r(M)}(M)=1$ as well. These cryptomorphisms constrain the structure of $\mathcal{L}(Q)\subseteq \mathcal{L}(M)$ considerably.
\begin{Cor} \label{Cor: monot_flats}
    Let $M$ and $Q$ be matroids on the same ground set $E$. If $M\twoheadrightarrow Q$, then $W_{k}(Q)\leq W_{k}(M)$ for all $0\leq k\leq r(M)$.
\end{Cor}
\begin{proof}
    By \Cref{Prop: quotient_crypto}, for every $X\subseteq E$, we have $r_Q(X)\leq r_M(X)$. For any $F\in \mathcal{L}_k(Q)$, we can choose a \emph{basis} $B_F$, which is an independent set in $Q$ satisfying $$r_Q(B_F)=|B_F|=k, \quad \text{cl}_Q(B_F)=F$$ Since $r_Q(X)\leq r_M(X)$, the following inequalities hold: $$k=r_Q(B_F)\leq r_M(B_F)\leq |B_F|=k$$ which implies $r_M(B_F)=k$. We now claim that for flats $F,G\in \mathcal{L}_k(Q)$ and choices of basis $B_F,B_G$, $$\text{cl}_M(B_F)=\text{cl}_{M}(B_G)\Longrightarrow F=G.$$
    Assume by contradiction that $F\neq G$ but $\text{cl}_M(B_F)=\text{cl}_{M}(B_G)=Y\subseteq E$. Then  $r_M(B_F)=r_M(\text{cl}_M(B_F))=r_M(Y)=k$, and this implies that $$k=r_Q(B_F)\leq r_Q(Y)\leq r_M(Y)=k.$$ In other words, $r_Q(Y)=k$. Identical arguments apply for $B_G$ and this implies that $\text{cl}_Q(B_F)=\text{cl}_Q(B_G)$, which is equivalent to $$F=\text{cl}_Q(Y)=G.$$ We have therefore shown that the map $F\mapsto \text{cl}_M(B_F)$ is an injection $\mathcal{L}_k(Q)\to \mathcal{L}_k(M)$, and hence $W_{k}(Q)\leq W_{k}(M)$ for all $0\leq k\leq r(M)$.
\end{proof}
We now move to some fundamental properties of quotients. We say a quotient $Q$ is of \emph{corank $k$ in $M$} if $r(Q)=r(M)-k$. A quotient of corank 1 is called an \emph{elementary quotient}. The set of elementary quotients of $M$ is denoted by $\mathcal{Q}^1(M)$. A \emph{linear subclass} of $M$ is a subset $\mathcal{H} \subseteq \mathcal{L}^1(M)$ that satisfies the
following property:
if $H_1,H_2\in \mathcal{H}$ cover $H_1\cap H_2\in \mathcal{L}^2(M)$, then all hyperplanes covering $H_1\cap H_2$ lie in $\mathcal{H}$. We call $\mathcal{H}=\mathcal{L}^1(M)$ the \emph{trivial} linear subclass. To each linear subclass $\mathcal{H}_{Q}$ we can additionally associate its \emph{modular cut} $\mathcal{M}_Q$, defined by $$\mathcal{M}_{Q}=\{F\in \mathcal{L}(M) \mid \text{ every hyperplane containing } F \text{ lies in } \mathcal{H}_{Q}\}$$ Foundational work of Crapo showed that elementary quotients, non-trivial linear subclasses, and modular cuts are all equivalent frameworks.

\begin{Prop} [\cite{crapo1964single}]\label{Prop: elt_quotients}
 The map $h:\mathcal{Q}^{1}(M)\to 2^{\mathcal{L}^1(M)}$ given by $Q\mapsto \mathcal{L}^1(Q)\cap \mathcal{L}^1(M)$ is a bijection from the elementary quotients of $M$ to the proper linear subclasses of $M$. Furthermore, the rank function of $Q$ is determined entirely by $\mathcal{M}_Q$: $$r_Q(X)=\begin{cases}
    r_M(X)-1 & \text{if } \text{cl}_M(X) \in \mathcal{M}_{Q}\\
    r_M(X) & \text{else}
\end{cases},\quad X\subseteq E(M)$$
\end{Prop}

The previous proposition gives us an understanding of elementary quotients $M\twoheadrightarrow Q$ in terms of \emph{high-rank} flats of $M$. This is especially useful because any quotient of $M$ can be written as a chain of elementary quotients.

\begin{Prop} [\cite{oxley2006matroid}] \label{Prop: Higgs_Lifts}
    A matroid $Q$ is a corank-$k$ quotient of $M$ if and only if there exists a sequence $$M=M_0\twoheadrightarrow M_1\twoheadrightarrow ...\twoheadrightarrow M_k=Q$$ such that $M_{i}\twoheadrightarrow M_{i+1}$ is an elementary quotient for all $0\leq i< k$.
\end{Prop}

\begin{rmk}
    \Cref{Prop: Higgs_Lifts} is equivalent to the statement that $\mathcal{Q}(M)$ is a poset graded by rank.
\end{rmk}

Switching gears, we also show that quotients behave naturally with respect to direct sum. This will be the crucial observation for disproving unimodality of quotient counts.
\begin{Lm} \label{Lm: prod_of_quotients}
    Let $M, N, Q$ be matroids satisfying $E(Q) = E(M )\sqcup E(N)$. Then $M\oplus N \twoheadrightarrow Q$  if and only if $M \twoheadrightarrow Q|_{E(M)}$ and $N \twoheadrightarrow Q|_{E(N)}$.
\end{Lm}
\begin{proof}
 Assume $M\oplus N \twoheadrightarrow Q$. To see that $M \twoheadrightarrow Q|_{E(M)}$ and $N \twoheadrightarrow Q|_{E(N)}$, it suffices to show that restriction preserves the quotient property. This is easy to see from the cryptomorphism for quotients via rank (c.f. \Cref{Prop: quotient_crypto}). Namely, the rank inequality still is satisfied when restricting to possible subsets $X\subseteq Y\subseteq E(M)$ or $X\subseteq Y\subseteq E(N)$.

 Now assume $M \twoheadrightarrow Q|_{E(M)}$ and $N \twoheadrightarrow Q|_{E(N)}$, and choose some $F\in \mathcal{L}(Q)$. Since restriction preserves flats, both $F\cap E(M)$ and $F\cap E(N)$ are flats of $Q|_{E(M)}$ and $Q|_{E(N)}$, respectively. This implies that $F\cap E(M)$ is a flat of $M$, and $F\cap E(N)$ is a flat of $N$. Since the flats of $M\oplus N$ are precisely the sets $F_M\sqcup F_N$ for $F_M\in \mathcal{L}(M)$, $F_N\in \mathcal{L}(N)$, $$F=(F\cap E(M))\sqcup (F\cap E(N))\in \mathcal{L}(M\oplus N)$$ and this implies that every flat of $Q$ is a flat of $M\oplus N$. In other words, $M\oplus N \twoheadrightarrow Q$.
\end{proof}

Lastly, each matroid $M$ has a particularly important quotient of rank $k$, called its \emph{$k$th truncation} $Tr_k(M)$. This matroid is defined by the rank function $$r_{Tr_k(M)}(X)=\min(k, r_M(X))$$ and its proper flats are precisely the flats of $M$ of rank less than $k$. We also denote the truncation operator $\text{Tr}(-)$ as the operator defined by $M\mapsto Tr_{r(M)-1}(M)$. 

\begin{Lm} \label{truncation is realizable over infinite extension}
    Let $M$ be realizable over an infinite field $\mathbb{F}$. Then $Tr_k(M)$ is also realizable over $\mathbb{F}$.
\end{Lm}
\begin{proof}
The result for $\text{Tr}(M)$ is shown by \cite[Theorem 8]{jurrius2012truncation}. By examining the rank function of $Tr_k(M)$, it is obvious that $$Tr_k(M)=\text{Tr}(\text{Tr}(\text{Tr}(...\text{Tr}(M))))$$ (where we apply the truncation operator $r(M)-k$ times) and hence we can repeatedly apply the result to show that $Tr_k(M)$ is also realizable over $\mathbb{F}$.
\end{proof}

We conclude this section with a few examples of matroids and their quotients.

\begin{Ex} \label{Ex: Uniform_quotients}
Let $Q$ be a quotient of $U_{2,n}$. Then:
    \begin{enumerate}
        \item [1.] If $Q$ has rank 0, it consists entirely of loops.
        \item [2.] If $Q$ has rank 1, then it has either no loops or precisely one loop.
        \item [3.] If $Q$ has rank 2, then $Q=U_{2,n}$.
    \end{enumerate}

    All but the rank 1 case follow from previous discussion. This case can easily be seen by considering the rank cryptomorphism for quotients. If there were at least two distinct loops $i,j$ and one non-loop $k$, choosing $X=\{i,j\}$, $Y=\{i,j,k\}$ gives a violation $$r_M(Y)-r_M(X)\geq r_Q(Y)-r_Q(X)\Longrightarrow 2-2\geq 1-0,$$ a contradiction. Conversely, any rank 1 matroid with at most one loop is clearly a quotient, since the loop is a proper flat of $U_{2,n}$. Thus there are precisely $n+2$ quotients of $U_{2,n}$ with rank $\leq 1$.
\end{Ex}

Given a poset $P$, its \emph{opposite} poset, $P^{\text{op}}$ is the poset with the same underlying set as $P$, with order relation reversed.
\begin{Df}[Modularity]
    A matroid $M$ is \emph{modular} if there exists an involutive poset isomorphism $\mathcal{L}(M)\simeq \mathcal{L}(M)^{\text{op}}$.
\end{Df}
\begin{rmk}
    This definition of modularity is slightly stronger than the definition commonly found in the literature. A matroid $M$ is generally called modular if for every pair of flats $X,Y$ in $M$, $r_M(X)+r_M(Y)=r_M(X\cup Y)+ r_M(X\cap Y)$. Our definition is the one given in \cite{wang2026tropical}, and is required for the statement of $\Cref{Thm: Modular_Symm}$ to hold.
\end{rmk} 

This modularity condition is particularly rigid, and modular matroids are difficult to come by as a result. Examples of modular matroids include the \emph{free matroid} $U_{n,n}$ and the \emph{projective geometries} $PG(r-1,q)$. The latter family will be especially important for our counterexamples to top-heaviness. Modular matroids are imbued with remarkable properties, particularly with respect to their quotients.

\begin{Thm} [\cite{wang2026tropical}]\label{Thm: Modular_Symm}
    If $M$ is modular, then $\mathcal{Q}(M)\simeq \mathcal{Q}(M)^{\text{op}}$. As such, $q_k(M)=q_{d-k}(M)$ for all $1\leq k \leq d$.
\end{Thm}

\section{Unimodality}
For the rest of the section, fix a matroid $M$ on ground set $E(M)$, and $U_{2,n}$ on $[n]$. We show that a violation of unimodality for $W_k(M)$ implies non-unimodality for $q_k(M\oplus U_{2,n})$, provided $n$ is sufficiently large. Recall that two integer-valued functions $f(n)$ and $g(n)$ satisfy $g(n)=\Theta(f(n))$ if there exist constants $c_1,c_2>0$ and some integer $n_0>0$ such that $c_1f(n)\leq g(n)\leq c_2f(n)$ for all $n>n_0$.

\begin{Thm} \label{Thm: Unimodality}
    Let M be a matroid realizable over some field, and let $1\leq k\leq r(M)+2$ be an integer. Then, asymptotically as $n$ grows, we have $q_k(M \oplus U_{2,n}) = \Theta(n^{W_{k-1}(M)})$.
\end{Thm} Directly applying this theorem to the matroid in \Cref{Thm: DLLWPeaks} and taking $n$ to be sufficiently large yields the following corollary:

\begin{Cor}
    For each positive integer $m$, there is a realizable matroid $M$ whose sequence $q_k(M)$ has at least $m$ peaks.
\end{Cor}
\begin{proof}
    By \Cref{Thm: DLLWPeaks}, there exists a realizable matroid with at least $m$ peaks. Directly applying \Cref{Thm: Unimodality} to this matroid yields a matroid with at least $m$ peaks. Since both $U_{2,n}$ and $M$ are assumed to be realizable, so is their direct sum. For every index such that $W_i(M)=W_{i+1}(M)$ it's not a priori clear that $q_{i+1}(M \oplus U_{2,n})\neq q_{i+2}(M \oplus U_{2,n})$, hence we lose the ``exactness" of the theorem in \Cref{Thm: DLLWPeaks}.
\end{proof}

The proof of \Cref{Thm: Unimodality} breaks into two parts: the lower bound and the upper bound. The upper bound holds without using the realizability condition on $M$, whereas the lower bound makes heavy use of this assumption.

\subsection{The upper bound.} Our arguments rely on one final definition, this time from the realm of polymatroids.

\begin{Df} [\cite{fujishige2005submodular}]
    Let $M$ be a matroid on ground set $E(M)$, and let $P=A_1\sqcup ... \sqcup A_n$ be a partition of $E(M)$ into non-empty sets, and let $E(P)$ be the set of blocks of $P$. Then the \emph{aggregation of $M$ along $P$} is the polymatroid $\mathcal{P}$ on $E(P)$ with rank function $$r_{\mathcal{P}}(X)=r_M\left(\bigcup_{A_i\in X}A_i\right),\quad X\subseteq E(P)$$
\end{Df} 

\begin{rmk}
    Aggregation has appeared across the literature under many different names. \cite{hall2004fork} refers to them as \emph{partitioned matroids}, \cite{helgason2006aspects} introduced them as special cases of \emph{hypermatroids} and \cite{berczi2026interaction} (ironically, for our purposes) refers to them as \emph{quotients of matroids by partitions}. We prefer the term ``aggregation" because it evokes the idea that the elements of a block are being combined in the resulting polymatroid.
\end{rmk}

Choose a non-empty set $S\subseteq E(M)$, and partition $E(M)$ by choosing one block for $S$, and all other elements as singletons. We denote the aggregation of $M$ along this partition as $\mathcal{P}_S(M)$, and one can essentially think of it as the polymatroid that results from replacing the set $S$ with a new ``generic" element $s$ in $M$. For any quotient $M\oplus U_{2,n}\twoheadrightarrow Q$, we can consider the polymatroid $\mathcal{P}_{[n]}(Q)$, which has fixed ground set $E(M)\sqcup \{s\}$ and rank \emph{at most} $r(M)+2$. The key idea of the proof is that we can separate the quotients $Q$ of $M\oplus U_{2,n}$ into finitely many groups based on their aggregation $\mathcal{P}_{[n]}(Q)$, and then bound 
the size of each group separately. We will break these arguments into cases based on the value of $0\leq r_Q([n])
\leq 2$. We first tackle the case that $r_Q([n])\leq 1$.

\begin{Lm} \label{Lm: rk_leq1_poly}
Let $A$ and $B$ be disjoint, non-empty finite sets and $s$ be some unrelated element. Fix a polymatroid $P$ on $A\sqcup \{s\}$ and a matroid $N$ on $B$. Let $\Omega_P$ denote the collection of matroids $Q$ on $A\sqcup B$ satisfying $Q|_{B}=N$ and $\mathcal{P}_{B}(Q) = P$. If $r(N)\leq 1$, then $|\Omega_P|\leq 1$.  
\end{Lm}
\begin{proof}
    If such a $Q$ exists, then $r(N)\leq 1$ implies $r_P(a)\leq 1$ for every $a\in A\cup \{s\}$. This means $P=\mathcal{P}_{B}(Q)$ is in fact a matroid, and $Q$ can be obtained from $P$ by adding in the loops of $N$ as loops of $Q$, and replacing $s$ with the non-loops of $N$ in parallel.
\end{proof}

In the case that $M\oplus U_{2,n}\twoheadrightarrow Q$ and $r_Q([n])=2$, combining \Cref{Prop: quotient_crypto} and \Cref{Lm: prod_of_quotients} implies that $Q|_{[n]}= U_{2,n}$. This case is handled by the following lemma.

\begin{Lm} \label{Lm: rk_2_poly}
Let $A$ and $B$ be disjoint, non-empty finite sets and let $s$ be some unrelated element. Fix a polymatroid $P$ on $A\sqcup \{s\}$ such that $P|_A$ is a matroid and $r_P(s)=2$. Let $\Omega_P$ denote the collection of matroids $Q$ on $A\sqcup B$ such that $Q|_{B}$ is uniform of rank $2$, and $\mathcal{P}_{B}(Q)=P$. Then an element $Q\in \Omega_P$ is picked out (in $\Omega_P$) by the function $H\mapsto \text{cl}_Q(H)\cap B$ defined on the set of flats $H$ of $P|_A$ of rank $r(P)-1$.
\end{Lm}
\begin{proof}
    Recall that a set $S\subseteq A\sqcup B$ is \emph{spanning} in $Q$ if $r_Q(S)=r(Q)$, and a matroid is determined by its collection of spanning sets (the bases of $Q$ are precisely the minimal spanning sets with respect to inclusion). 

    Choose a set $S\subset E(Q)$. We will show that the function $$H\mapsto \text{cl}_Q(H)\cap B$$ contains enough information to determine if $S$ is spanning. In many cases, we won't even need this function to conclude that our set is spanning. We first observe that without loss of generality, $S\cap A$ can be assumed to be a flat of $Q|_A$. This is because $S$ is spanning in $Q$ if and only if the same is true of $S\cup \text{cl}_{P|_A}(S\cap A)$, since $\text{cl}_{P|_A}(S\cap A)\subseteq \text{cl}_Q(S)$. Now that we have constrained the structure of $H=S\cap A$, we move on to $S\cap B$. This breaks into three cases based on the size of this intersection. First, if $S\cap B =\emptyset$, then $S$ is spanning in $Q$ if and only if $r_{P|_A}(S)=r(P)$. If $|S \cap B| \geq 2$, then $S$ is spanning in $Q$ if and only if $(S\setminus B) \cup \{s\}$ is spanning in P (meaning it has rank $r(P)$). Now suppose $|S\cap B|=1$. Then $S=H\cup \{b\}$ for some $b \in B$ and some flat $H$ of $P|_A = Q|_A$. If $\text{rk}_{P|_A}(H) = r(P)$, then $S$ is spanning in $Q$. Now assume that $\text{rk}_{P|_A}(H) \leq r(P)-2$. This would imply that $$\text{rk}_{Q}(S) \leq \text{rk}_{P|_A}(H)+1<r(P)=r(Q)$$ so $S$ clearly cannot be spanning. This leaves the case that $\text{rk}_{P|_A}(H) = r(P)-1$. This is where the function $H\mapsto \text{cl}_Q(H)\cap B$ enters; whether or not $S$ is spanning depends entirely on the ``location" of $b\in B$. If $b\in \text{cl}_Q(H) \cap B$, then $$\text{rk}_Q(S) = \text{rk}_Q(H) = \text{rk}_{P|_A}(H) < r(P)$$ and $S$ is not spanning. If $b\notin \text{cl}_Q(H) \cap B$, then $$\text{rk}_Q(S) = \text{rk}_Q(H)+1 = r(P)$$ and $S$ is spanning. This means we have determined precisely the spanning sets of $Q$, and hence $Q$ is uniquely specified by the function $H\mapsto \text{cl}_Q(H)\cap B$.
\end{proof}

Armed with these results, we can prove the upper bound in \Cref{Thm: Unimodality}. 
\begin{Prop} \label{Prop: Upper bound}
    Let M be a matroid and let $1\leq k\leq r(M)+2$ be an integer. Then there exists some $C>0$ such that $q_k(M \oplus U_{2,n})\leq C\cdot n^{W_{k-1}(M)}$.
\end{Prop}
\begin{proof}
Now fix a polymatroid $P$ on $E(M)\sqcup \{s\}$ of rank $k$. We will first bound the number of quotients $M\oplus U_{2,n}\twoheadrightarrow Q$ such that $\mathcal{P}_{[n]}(Q)=P$. By \Cref{Lm: prod_of_quotients}, any such $Q$ must satisfy $M\twoheadrightarrow Q|_{E(M)}$ and $U_{2,n}\twoheadrightarrow Q|_{[n]}$ and $M\twoheadrightarrow P|_{E(M)}$. As above, our proof breaks into cases based on the value of $r_Q([n])$. 

If $r_Q([n])\leq 1$, choose $Q|_{[n]}$ to be a fixed (rank 1) quotient $V$ of $U_{2,n}$. Then \Cref{Lm: rk_leq1_poly} implies there can exist at most one matroid $Q$ with $\mathcal{P}_{[n]}(Q)=P$ and $Q|_{[n]}=V$. Since $U_{2,n}$ has $n+2$ different quotients of rank $\leq 1$ (c.f. \Cref{Ex: Uniform_quotients}), there can be at most $n+2$ quotients arising from this case.

If $r_Q([n])=2$, then \Cref{Prop: quotient_crypto} implies that $Q|_{[n]}= U_{2,n}$. We can apply \Cref{Lm: rk_2_poly} to see that the matroid $Q$ is determined by the function $\omega: \mathcal{L}_{k-1}(P|_{E(M)})\to \mathcal{L}(U_{2,n})$ given by $$H \mapsto \text{cl}_Q(H)\cap [n].$$ The number of quotients $Q$ with $\mathcal{P}_{[n]}(Q)=P$ is therefore bounded above by the number of possible functions $\omega$. Since $|\mathcal{L}_{k-1}(P|_{E(M)})|=W_{k-1}(P|_{E(M)})$ and $|\mathcal{L}(U_{2,n})|=n+2$, this becomes $$(n+2)^{W_{k-1}(P|_{E(M)})}\leq (n+2)^{W_{k-1}(M)}=O(n^{W_{k-1}(M)}),$$ where the first inequality follows because of \Cref{Cor: monot_flats}, which implies that $W_{k-1}(P|_{E(M)})\leq W_{k-1}(M)$ for all $1\leq k \leq r(M)+1$.

Finally, we observe that the number of integral polymatroids of a fixed rank and ground set (neither of which depend on $n$) is finite, so we can sum over all such polymatroids without affecting our final asymptotic bound. The $n+2$ quotients contributed by rank 1 case are negligible outside of the cases that $k$ is extremal ($k=1,r(M)+1,r(M)+2$). When $k=r(M)+2$, there are no flats of $M$ with rank $k-1$, and there is a unique quotient. Our asymptotic estimate still holds in this case. When $k=1, r(M)+1$, there is a unique flat in $M$ of rank $k-1$, and when $W_{k-1}(M)\geq 1$ the asymptotic estimate still holds in this case. In all other cases, we have shown that $q_k(M\oplus U_{2,n})=O(n^{W_{k-1}(M)})$, as desired.
\end{proof}

\subsection{Lower bound}
The lower bound is comparatively simple, as it leverages the realizability of $M$ over some field. Our main statement hinges on a statement about hyperplanes:

\begin{Prop} \label{There exists 2-dimensional subspace}
Let $H_1,H_2,...,H_k$ be distinct hyperplanes in $K^r$, where $k,r\geq 2$ and $K$ is an infinite field. Then there exists a two-dimensional subspace $L\subset K^r$ such that $L\cap H_i$ are pairwise distinct one-dimensional subspaces.
\end{Prop}
\begin{proof}
For each $H_i$ there exists a linear form $\ell_i: K^r\to K$ such that $H_i=\ker(\ell_i)$. For $i \neq j$, the forms $\ell_i, \ell_j$ are clearly
linearly independent, otherwise their kernels would agree. In order to construct $L$, we view it as the span of two non-zero vectors $u,v\in K^r$. We first choose $u\in K^r$ such that $\ell_i(u)\neq 0$ for all $1\leq i\leq k$. This necessarily exists because $K$ is infinite and there are only finitely many $H_i$ of co-dimension 1 to avoid. Now choose $v\in K^r$ so that $$\ell_i(u)\ell_j(v)-\ell_i(v)\ell_j(u)\neq 0$$ for all $i\neq j$. Again, this is a finite set of hyperplanes over an infinite field, so such a $v$ necessarily exists. Furthermore, $u$ and $v$ are linearly independent (otherwise one of the above polynomials would vanish), and on
$L= \text{span}_K\{u,v\}$  the restrictions of every pair of forms are independent. Each restriction $L\cap H_i=\text{ker}(\ell_i|_L)$
is one dimensional, because $\ell_i(u)\neq 0$, and by our choice of $v$, $L\cap H_i$ and $L\cap H_j$ are distinct for all $i\neq j$.
\end{proof}

\begin{Prop}
    Let M be a realizable matroid and let $1\leq k \leq r(M)+2$ be an integer. Then $\frac{1}{2}n^{W_{k-1}(M)}\leq q_{k}(M\oplus U_{2,n})$ for sufficiently large $n$.
\end{Prop}
\begin{proof}
    We first tackle the extremal cases, where $k=1,r(M)+1,r(M)+2$. When $k=r(M)+2$, there is no flat of rank $k-1$ in $M$, but one rank $k$ quotient of $M\oplus U_{2,n}$. The lower bound immediately holds. When $k=1$ (resp. $r(M)+1$), take the rank 0 (resp. full rank) quotient of $M$ and choose a rank 1 quotient $V$ of $U_{2,n}$. This yields a quotient of $M\oplus U_{2,n}$. There are $n+1$ such choices for $V$, so this process yields $n+1$ distinct quotients. Since there is a unique flat of rank $0$ and $r(M)$ in $M$, this yields the desired conclusion in all cases.
    
    Now assume $2\leq k \leq r(M)$. In fact, we show there are at least $\frac{1}{2}n^{W_{k-1}(M)}$ \emph{realizable} rank $k$ quotients $M\oplus U_{2,n}\twoheadrightarrow Q$ which satisfy $Q|_{[n]}=U_{2,n}$ and $Q|_{E(M)}=Tr_k(M)$. Since any realizable matroid admits a realization over an infinite field, we can assume that $M$ is realizable over some infinite field $K$. \Cref{truncation is realizable over infinite extension} then implies that $Tr_k(M)$ is also realizable over $K$. Choose vectors $v_1,...,v_m\in K^k$ realizing $Tr_k(M)$. For any flat $F\in \mathcal{L}_{k-1}(M)$, define the hyperplane $$H_F=\text{span}_K(v_f \mid f\in F).$$ We note that if $F_1, F_2\in \mathcal{L}_{k-1}(M)$ are distinct, then $r_M(F_1\cup F_2)\geq k$ and therefore $\text{dim}(\text{span}_K(v_f \mid f\in F_1\cup F_2))=k$. This implies that $H_{F_1}\neq H_{F_2}$, and hence we have constructed $W_{k-1}(M)$ distinct hyperplanes of $Tr_k(M)$ in $K^k$. There are also at least two distinct flats of rank $k-1$ in $Tr_k(M)$, as we can pick a basis of $M$ and restrict to two distinct subsets of size $k-1$. We can then apply \Cref{There exists 2-dimensional subspace} (since we assumed that $2\leq k\leq r$) to see that there exists a 2-dimensional subspace $L\subseteq K^k$ with $L\cap H_{F_1}\neq L\cap H_{F_2}$ for any $F_1,F_2\in \mathcal{L}_{k-1}(M)$. We consider $L$ as a projective line, and view $p_F=L\cap H_{F}$ as a collection of pairwise distinct projective points on $L$. 

    Now suppose that $n\geq W_{k-1}(M)$. For every injection $\phi: \mathcal{L}_{k-1}(M)\hookrightarrow [n]$, we label the point $p_F$ by $\phi(F)$. Place the other labels of $[n]$ at pairwise distinct projective points of $L\setminus \{p_F \mid F \in\mathcal{L}_{k-1}(M)\}$. We are able to do this without trouble because $K$ is infinite. For every point $i\in [n]\setminus \phi(\mathcal{L}_{k-1}(M))$, choose a non-zero vector $v_i$ realizing $p_i$. Taking the union of these vectors with $\{v_1,...,v_m\}$ gives us a new (rank $k$) matroid $Q_\phi$, depending on our choice of injection. By construction, $Q_\phi|_{[n]}=U_{2,n}$ (since each projective point was pairwise distinct), and $Q_\phi|_{E(M)}=Tr_k(M)$. \Cref{Lm: prod_of_quotients} therefore implies that $M\oplus U_{2,n}\twoheadrightarrow Q_\phi$.

    It remains to show that different choices of $\phi$ yield different quotients of $M\oplus U_{2,n}$. This is easy to see from the observation that for any $F\in \mathcal{L}_{k-1}(M)$ and $i\in [n]$, $$\phi(F)=i\Longleftrightarrow r_{Q_\phi}(F\cup i)=k-1.$$ This holds because $v_i\in H_F$ precisely when its projective point is $p_F$. As such, changing $\phi$ changes the rank function, and we have shown that different injections yield different rank $k$ quotients of $M\oplus U_{2,n}$.

    The number of injections $\mathcal{L}_{k-1}(M)\hookrightarrow [n]$ is precisely the falling factorial $(n)_{W_{k-1}(M)}$, and we therefore have shown that $$q_k(M\oplus U_{2,n})\geq (n)_{W_{k-1}(M)}$$ To turn this into our desired lower bound, observe that the following inequalities hold for $n\geq D(D-1)$, $$(n)_D= n^{D}\prod_{i=0}^{D-1}(1-\frac{i}{n})\geq n^{D}(1-\sum_{i=0}^{D-1}\frac{i}{n})\geq n^{D}(1-\frac{D(D-1)}{2n})\geq \frac{1}{2}n^D.$$ Plugging in for $D=W_{k-1}(M)$ implies that $$\frac{1}{2}n^{W_{k-1}(M)}\leq q_k(M\oplus U_{2,n})$$ for $n$ sufficiently large.
\end{proof}

\section{Top-heaviness}
In contrast to the previous section, we derive a small, explicit example where top-heaviness fails. 

\begin{Thm} \label{Thm: Top_Heaviness}
    Let $q$ be a prime power and $d\geq 5$. Then for every point $p$ of $PG(d-1,q)$, the
matroid $M=PG(d-1,q)\setminus p$ satisfies $q_2(M)> q_{d-2}(M)$.
\end{Thm}

The idea behind this counterexample is simple; since $G=PG(d-1,q)$ is modular, \Cref{Thm: Modular_Symm} implies that $q_k(G)=q_{d-k}(G)$ for all $k$. The fact that equality holds means that $G$ is already quite close to violating top-heaviness. The naive hope is that ``perturbing" $G$ in the right way could cause a violation. It turns out that deleting a point from $G$ doesn't affect quotient counts for any rank $\geq 3$, but it increases the number of rank 2 quotients. For a matroid $M$, we define the subposets $$\mathcal{Q}^{\leq k}(M)=\{Q\in \mathcal{Q}(M) \mid r(M)-r(Q)\leq k\}$$ $$\mathcal{L}^{\leq k}(M)=\{F\in \mathcal{L}(M) \mid r(M)-r_M(F)\leq k\}$$ to be the subposet of quotients (resp. flats) of corank $\leq k$. Before proving our main result, we prove the following auxiliary statement:

\begin{Lm} \label{Lm: desc_quotients}
    Let $M, M'$ be matroids, fix an integer $k\geq 1$, and let $\phi: \mathcal{L}^{\leq k+1}(M) \to \mathcal{L}^{\leq k+1}(M')$ be a poset isomorphism. Then for every elementary
quotient $M \twoheadrightarrow Q$, there exists a unique elementary quotient $M' \twoheadrightarrow Q'$ such that $\phi$ restricts to an isomorphism $\mathcal{L}^{\leq k}(Q) \to \mathcal{L}^{\leq k}(Q')$.
\end{Lm}
\begin{proof}
Let $Q$ be an elementary quotient of $M$. By \Cref{Prop: elt_quotients}, $Q$ is determined by an associated proper linear subclass $\mathcal{H}_{Q}$, and an associated modular cut $\mathcal{M}_{Q}$. Since linear subclasses are entirely determined by the flats of corank $\leq 2$ in $M$ (and we assumed $k\geq 1$), $\mathcal{H}_{Q}$ is determined by $\mathcal{L}^{\leq k+1}(M)$. Furthermore, the corank of any flat $F\in \mathcal{L}^{\leq k+1}(M)$ is determined by the length of its upper interval $[F,M]$, and since $\phi$ is a poset isomorphism, it must preserve that length. Thus $\phi$ is corank-preserving. Since $\phi$ also preserves containment relations, it sends $\mathcal{H}_Q$ to a proper linear subclass $\mathcal{H}\subseteq \mathcal{L}^{\leq k+1}(M')$. By definition of linear subclasses, $\mathcal{H}=\mathcal{H}_{Q'}$ for a unique elementary quotient $M'\twoheadrightarrow Q'$. It remains to show that $\phi|_{\mathcal{L}^{\leq k}(Q)}:\mathcal{L}^{\leq k}(Q) \to \mathcal{L}^{\leq k}(Q')$ is an isomorphism.

At the moment it's not even clear that $\mathcal{L}^{\leq k}(Q)\subseteq\mathcal{L}^{\leq k+1}(M)$. As discussed above, an elementary quotient $Q$ is also determined by its modular cut $\mathcal{M}_{Q}$: $$\mathcal{M}_{Q}=\{F\in \mathcal{L}(M) \mid \text{ every hyperplane containing } F \text{ lies in } \mathcal{H}_{Q}\}$$ Observe that whether or not $F$ lies in $\mathcal{M}_{Q}$ is also entirely determined by the interval $[F,M]$. Most importantly, the rank function of $Q$ is determined by $\mathcal{M}_Q$: 
$$r_Q(X)=\begin{cases}
    r_M(X)-1 & \text{if } \text{cl}_M(X) \in \mathcal{M}_{Q}\\
    r_M(X) & \text{else}
\end{cases},\quad X\subseteq E(M)$$

Since $Q$ is an elementary quotient, we can substitute in $r(Q)=r(M)-1$ to get an expression for the corank function of $Q$:  
$$\text{crk}_Q(X)=\begin{cases}
    \text{crk}_M(X) & \text{if } \text{cl}_M(X) \in \mathcal{M}_{Q}\\
    \text{crk}_M(X)-1 & \text{else}
\end{cases},\quad X\subseteq E(M)$$

Since $\text{crk}_Q(X)\leq k$ for every flat of $\mathcal{L}^{\leq k}(Q)$, this gives us the desired inclusion $\mathcal{L}^{\leq k}(Q)\subseteq\mathcal{L}^{\leq k+1}(M)$. It remains to show that $\phi$ sends the flats of $Q$ to the flats of $Q'$. By studying the above rank function for $Q$, we can extract the lattice of flats $\mathcal{L}(Q)$ in terms of $\mathcal{M}_Q$ as well:
$$F\in \mathcal{L}(Q) \Longleftrightarrow F\in \mathcal{M}_{Q} \text{ or no flat covering } F \text{ belongs to } \mathcal{M}_Q.$$

We observe that the property determining when $F\in \mathcal{M}_Q$ can be expressed purely in the information given in $\mathcal{L}^{\leq k+1}(M)$. This is true because $\mathcal{M}_Q$ is upwards-closed, and the interval $[F,M]$ entirely determines when $F\in \mathcal{M}_Q$. To summarize, we have shown that $\mathcal{L}^{\leq k}(Q)$ and the corank of those flats can be entirely determined by $\mathcal{L}^{\leq k+1}(M)$. Since $\phi$ preserves order, coranks, and hyperplanes, it sends $\mathcal{H}_Q$ to $\mathcal{H}_{Q'}$, and this implies that $$F\in \mathcal{L}^{\leq k}(Q) \Longleftrightarrow \phi(F)\in \mathcal{L}^{\leq k}(Q').$$
\end{proof}

Inducting on the previous statement yields a surprising relationship between the upper parts of $\mathcal{Q}(M)$ and $\mathcal{L}(M)$.

\begin{Prop} \label{Prop: upper_flats_quotients}
    Let $M,M'$ be matroids, and fix an integer $k\geq 0$. Then any isomorphism $\phi: \mathcal{L}^{\leq k+1}(M) \to \mathcal{L}^{\leq k+1}(M')$ induces a (corank-preserving) poset isomorphism
$$\varphi: \mathcal{Q}^{\leq k}(M) \to \mathcal{Q}^{\leq k}(M').$$ Furthermore, $\varphi(Q)$ is determined by the image of $\mathcal{L}^1(\varphi(Q)) = \varphi(\mathcal{L}^1(Q))$.
\end{Prop}
\begin{proof}
    Our first goal will be to show that $\varphi$ is a well-defined bijection. We then will show that it respects the partial order on quotients. 
    
    Choose an arbitrary quotient $Q\in \mathcal{Q}^{\leq k}(M)$ of corank $d$. Then by \Cref{Prop: Higgs_Lifts}, there exists a sequence of elementary quotients $M=M_0\twoheadrightarrow M_1\twoheadrightarrow... \twoheadrightarrow M_d = Q$. Now assume there exists a matroid $M_0'=M'$ and an accompanying poset isomorphism $\phi_0=\phi: \mathcal{L}^{\leq k+1}(M_0) \to \mathcal{L}^{\leq k+1}(M'_0).$ Iteratively applying \Cref{Lm: desc_quotients} $i$ times shows there exists a quotient $M'\twoheadrightarrow M_i'$ and an isomorphism $$\phi_i: \mathcal{L}^{\leq k+1-i}(M_i) \to \mathcal{L}^{\leq k+1-i}(M'_i).$$ This process terminates at some quotient $M_d'=Q'$ and an isomorphism $$\phi_d: \mathcal{L}^{\leq k+1-d}(Q) \to \mathcal{L}^{\leq k+1-d}(Q').$$  
    This is the case because $i\leq d\leq k$, so $k+1-d\geq 1$ allows us to apply the lemma at this final step. Since $\phi_d$ is a corank-preserving restriction of $\phi$ to the flats $\mathcal{L}^{\leq k+1-d}(Q)$, and $k+1-d\geq 1$ again implies that  $\mathcal{L}^{\leq 1}(Q)\subseteq \mathcal{L}^{\leq k+1-d}(Q)$, the hyperplanes of $Q'$ are uniquely determined to be $\mathcal{L}^1(Q')=\phi(\mathcal{L}^1(Q))$. Since a matroid on a fixed ground set is determined by its hyperplanes, $Q'$ only depends on $Q$ and $\phi$, and not our factorization of $M\twoheadrightarrow Q$ into elementary quotients. Define $\varphi(Q)=Q'$. By applying the same logic to $\phi^{-1}$, the hyperplanes of $Q$ map to $\mathcal{L}^1(Q)=\phi^{-1}(\mathcal{L}^1(Q'))$ and these likewise determine the matroid $Q$. We define $\varphi^{-1}(Q')=Q$, and composing these two morphisms (in either order) gives the identity. As such, $\varphi$ is a bijection. It also necessarily preserves corank, as applying the lemma $d$ times constructs a factorization of $M'\twoheadrightarrow Q'$ into $d$ elementary quotients. We have $\varphi(M)=M'$ as a result.

    It remains to show that for any quotient $P$ with $Q\leq P$, then $\varphi(Q)\leq \varphi(P)$. Assume $Q\leq P$. Then there exist factorizations of $M\twoheadrightarrow P$ and $P\twoheadrightarrow Q$, and concatenating them yields a factorization $$M\twoheadrightarrow...\twoheadrightarrow P \twoheadrightarrow...\twoheadrightarrow Q$$ into elementary quotients. We already showed that $\varphi(Q)$ is independent of the factorization of $M\twoheadrightarrow Q$, and the same is true of $\varphi(P)$ (since they both have corank $\leq k$). Iteratively applying \Cref{Lm: desc_quotients} therefore creates a chain of elementary quotients $$M'\twoheadrightarrow ...\twoheadrightarrow\varphi(P)\twoheadrightarrow ...\twoheadrightarrow\varphi(Q)$$ and hence $\varphi(Q)\leq \varphi(P)$. Identical arguments apply for $\varphi^{-1}$ so we have shown it is order preserving, and thus $\varphi$ is a corank-preserving poset isomorphism.
\end{proof}

We now have the background required to prove \Cref{Thm: Top_Heaviness}.

\begin{proof}
 Let $G=PG(d-1,q)$ and $M=G\setminus p$. Since $G$ is modular, \Cref{Thm: Modular_Symm} implies that $q_k(G)=q_{d-k}(G)$ for every $0\leq k \leq d$. We now observe that every flat $F\in \mathcal{L}^{\leq d-2}(G)$ is spanned by the elements of $F\setminus p$, regardless of our choice of $p$. This yields an isomorphism $\mathcal{L}^{\leq d-2}(M)\simeq \mathcal{L}^{\leq d-2}(G)$. \Cref{Prop: upper_flats_quotients} then constructs a corank-preserving isomorphism $$\mathcal{Q}^{\leq d-3}(M)\simeq \mathcal{Q}^{\leq d-3}(G),$$ Combining these observations gives the equalities $q_k(M)=q_k(G)=q_{d-k}(G)$ for all $3\leq k \leq d$. Since $d\geq 5$, we have $d-2\geq 3$, and this implies that $$q_{d-2}(M)=q_{d-2}(G)= q_2(G).$$ 
    
    To prove the violation of top-heaviness, it therefore suffices to show that $q_2(M)>q_2(G)$. We first show that $q_2(M)\geq q_2(G)$. This follows from showing that the map $f: \mathcal{Q}_2(G)\to \mathcal{Q}_2(M)$ given by $Q\mapsto Q|_{E(G)\setminus p}$ is a well-defined injection. Clearly, $M\twoheadrightarrow Q|_{E(G)\setminus p}$ (restricting the ground set does not affect the rank cryptomorphism defining quotients in \Cref{Prop: quotient_crypto}). We also claim that $r_M(Q|_{E(G)\setminus p})=2$; if it were rank 1, this would imply that $p$ is a coloop and therefore that $E(G)\setminus p$ were a hyperplane in $Q$. But this would imply that $E(G)\setminus p$ were a flat of $G$, which clearly cannot be the case since deleting a single element fixes the rank of a projective geometry. Thus we have shown that $f$ is well defined. 
    
    Let $R=f(Q)$ have rank $2$. We now show that the rank function of $R$ determines the rank function of $Q$. Let $X\subseteq E(G)\setminus p$ be non-empty. Clearly, we can understand $r_Q(X)$ in terms of $R$. Now choose some $x\in X$, and consider the line $\ell$ through the points $x,p$ in $G$. Since $\ell \setminus p$ spans $\ell$ and it consists of multiple points, we can conclude that $X\cup \{p\}$ and $X\cup (\ell\setminus \{p\})$ span the same flat in $G$. This implies (by \Cref{Prop: quotient_crypto}) that $$X\cup \{p\}\subseteq \text{cl}_{G}(X\cup (\ell\setminus \{p\}))\subseteq \text{cl}_{Q}(X\cup (\ell\setminus \{p\})), \quad X\cup (\ell\setminus \{p\})\subseteq \text{cl}_{G}(X\cup \{p\})\subseteq \text{cl}_{Q}(X\cup \{p\})$$ or in other words, $X\cup \{p\}$ and $X\cup (\ell\setminus \{p\})$ have the same closure in $Q$. This allows us to determine $r_Q(X\cup \{p\})$ from $R$ alone. Lastly, when $X$ is empty, $p$ is not a loop precisely when there exists a line $\ell$ through $p$ of rank 2 in $R$. This is the case because $\ell$ is spanned by $p$ and any other point on it, and conversely $p$ belongs to a basis of $Q$ as a non-loop. In this way, the rank function of $Q$ is entirely determined by the rank function of $R$. As such, our map $f$ is an injection, and we have shown that $$q_2(M)\geq q_2(G)=q_{d-2}(M).$$
    
    It remains to construct a rank 2 quotient $M\twoheadrightarrow Q$ which does not extend to $G$. For a line $\ell\in \mathcal{L}_2(G)$ containing $p$, we define $D_{\ell}=\ell\setminus p$. Consider two distinct lines $\ell_1,\ell_2$ containing $p$. Since the flats of $M=G\setminus p$ are the sets of the form $F\setminus p$ for $F\in \mathcal{L}(G)$, both $D_{\ell_1}$ and $D_{\ell_2}$ are flats of $M$. Now define a rank 2 matroid $Q$ by assigning all elements of $D_{\ell_1}$ (resp. $D_{\ell_2}$) to form classes of parallel elements, and assigning all other elements to be singleton (non-loop) classes. The flats of $Q$ are all flats of $M$ by the previous discussion, so $M\twoheadrightarrow Q$. 
    
    Assume by contradiction that there existed some $G\twoheadrightarrow Q'$ with $Q'|_{E(G)\setminus p}=Q$. Since $D_{\ell_1}$ and $D_{\ell_2}$ are both parallel classes in $Q$, this implies $$r_{Q'}(D_{\ell_1})=r_{Q'}(D_{\ell_2})=1$$ Since both $D_{\ell_1}$ and $D_{\ell_2}$ span $p$ in $G$, this implies that $r_{Q'}(\ell_1)=r_{Q'}(\ell_2)=1$ and since $p$ is spanned by two distinct parallel classes, $p$ must be a loop of $Q'$. Now choose a third line $\ell_3$ through $p$ which is distinct from $\ell_1,\ell_2$. Then there exist two points $x,y\in \ell_3\setminus p$ which form singleton classes in the quotient $Q$. As such, $r_Q(\{x,y\})=2$. But over $G$, we have $$y\in \text{cl}_G(\{x,p\})\subseteq \text{cl}_{Q'}(\{x,p\})=\text{cl}_{Q'}(\{x\}).$$
    This gives $r_Q(\{x,y\})=r_{Q'}(\{x,y\})=1$, a contradiction, since we assumed $y$ was not a loop in $Q$. We have therefore shown a strict inequality $$q_{2}(M)> q_{d-2}(M).$$
\end{proof}

\printbibliography

\end{document}